\documentclass[a4paper,10pt]{amsart}
\usepackage{CJK}
\usepackage{amsmath}
\usepackage{amsthm,amssymb,latexsym}
\usepackage{amsmath,calligra}
\usepackage{url,ae}
\usepackage{t1enc}
\usepackage{eufrak}
\usepackage{amsfonts}
\usepackage{graphicx}
\usepackage{hyperref}
\usepackage{amscd}
\usepackage{mathrsfs}
\usepackage{bbm}
\usepackage{color}
\usepackage{txfonts}

\usepackage{multirow,makecell}
\usepackage[format=hang,font=small,textfont=it]{caption}

\newtheorem{theorem}{Theorem}[section]

\newtheorem{proposition}[theorem]{Proposition}

\newtheorem{question}{Question}[section]

\theoremstyle{definition}
\newtheorem{definition}[theorem]{Definition}
\theoremstyle{remark}
\newtheorem{remark}[theorem]{Remark}
\theoremstyle{example}

\newcommand{\be}{\begin{equation}}
\newcommand{\ee}{\end{equation}}
\newcommand{\bea}{\begin{eqnarray}}
\newcommand{\eea}{\end{eqnarray}}
\newcommand{\ben}{\begin{eqnarray*}}
\newcommand{\een}{\end{eqnarray*}}
\newcommand{\bet}{\begin{equation}
\begin{split}}
\newcommand{\eet}{\end{split}
\end{equation}}

\DeclareMathOperator{\Ext}{\mathscr{E}\text{\kern -3pt {\calligra\Large xt}}\,\,}

\begin{document}

\title[On a question of Koll\'{a}r]
      {On a question of Koll\'{a}r}

\author{Zhenqian Li}

\date{\today}
\thanks{E-mail: lizhenqian@amss.ac.cn}

\begin{abstract}
In this note, we establish a generalized analytic inversion of adjunction via the Nadel-Ohsawa multiplier/adjoint ideal sheaves associated to plurisubharmonic (psh) functions for log pairs, by which we answer a question of Koll\'{a}r in full generality.
\end{abstract}

\maketitle

\section{Introduction}

Throughout this note, all complex spaces are always assumed to be reduced and paracompact unless otherwise mentioned; we mainly refer to \cite{GR84, Richberg68} for basic references on the theory of complex spaces. For various terminologies and basic facts in algebraic geometry, we refer to the standard expositions \cite{Kollar97, Kollar_MMP, K-M98, La04}.

Let $X$ be an $n$-dimensional normal complex space with a canonical (Weil) divisor/class $K_X$ and $\Delta=\sum d_k\Delta_k$ a Weil $\mathbb{R}$-divisor on $X$, where $\Delta_k$ are distinct prime divisors. The pair $(X,\Delta)$ is called to be a \emph{log pair} if the $\mathbb{R}$-divisor $K_X+\Delta$ is an $\mathbb{R}$-Cartier divisor on $X$. Let $\pi:\widetilde X\to X$ be a log resolution of $(X,\Delta)$, i.e., a proper modification $\pi:\widetilde X\to X$ such that $\widetilde X$ is smooth, the exceptional locus $\text{Ex}(\pi)$ of $\pi$ is a divisor and $\pi_*^{-1}\Delta+\text{Ex}(\pi)$ has simple normal crossings (SNC). Write $K_{\widetilde X}=\pi^*(K_X+\Delta)+\sum_{i}a(E_i,X,\Delta)E_i$, where $a(E_i,X,\Delta)\in\mathbb{R}$ and $E_i\subseteq\widetilde X$ are distinct prime divisors. We say that a log pair $(X,\Delta)$ is \emph{KLT} (respectively, \emph{LC}) if $a(E_i,X,\Delta)>-1$ (respectively, $a(E_i,X,\Delta)\geq-1$) for all $i$. Similarly, a log pair $(X,\Delta)$ is called to be \emph{PLT} if $d_k\leq1$ and for all log resolutions $\pi:\widetilde X\to X$ of $(X,\Delta)$ we have $a(E_i,X,\Delta)>-1$ for every exceptional divisor $E_i$.

Suppose that $(X,S+B)$ is a log pair, where $S\subseteq X$ is a reduced Weil divisor (complex subspace of pure codimension one) with $\xi:\widehat S\to S$ the normalization of $S$, and $B\subseteq X$ is an $\mathbb{R}$-divisor on $X$ which has no irreducible components in common with $S$. Then, there exists a unique $\mathbb{R}$-divisor $\text{Diff}_{\widehat S}(B)$ (called the \emph{different} of $B$) on $\widehat S$ such that $(K_X+S+B)|_{\widehat S}=K_{\widehat S}+\text{Diff}_{\widehat S}(B)$. In addition, if $\widetilde S=\pi_*^{-1}S$ is the strict transform (may be disconnected) of $S$ on $\widetilde X$ and $\mu:\widetilde S\to\widehat S$ is the induced morphism, then $\text{Diff}_{\widehat S}(B)=\mu_*\big((\Delta_{\widetilde X}-\widetilde S)|_{\widetilde S}\big)$, where $\Delta_{\widetilde X}=\pi^*(K_X+S+B)-K_{\widetilde X}$; see \cite{Kollar92, Kollar_MMP} for more details.

\subsection{Multiplier and adjoint ideal sheaves for log pairs}

Let $(X,S+B)$ be a log pair as above and $\varphi\in L_\text{loc}^1(X_\text{reg})$ be a weight function on $X$ such that $\varphi|_S\not\equiv-\infty$ on every irreducible component of $S$. In the present note, by introducing an Ohsawa-type measure on $S$ via adapted volume forms for the pair $(X,S+B)$, we define the Nadel-Ohsawa multiplier ideal sheaves on $S$ associated to the weight $\varphi$ for the pair $(X,S+B)$; see Definition \ref{MIS_NO} for more details. Relying on the notion of Nadel-Ohsawa multiplier ideal sheaf and an $L^2$ extension theorem of measure version, we are able to obtain a reasonable generalization of adjoint ideal sheaves to the analytic setting for log pairs and establish the adjunction exact sequence as follows:

\begin{theorem} \label{Aadjoint}
With the same notations as above, $B$ an effective $\mathbb{R}$-divisor and $\varphi\in\emph{QPsh}(X)$ a quasi-psh function on $X$. Then, there exists an ideal sheaf $$\emph{Adj}_S(X,S+B;\varphi)\subseteq\mathcal{O}_X,$$
called the \emph{analytic adjoint ideal sheaf} associated to the triple $(X,S+B;\varphi)$ along $S$, sitting in an exact sequence:
$$0\longrightarrow\mathscr{I}(X,S+B;\varphi)\stackrel{\iota}{\longrightarrow}\emph{Adj}_S(X,S+B;\varphi)\stackrel{\rho}{\longrightarrow} i_*\mathscr{I}_\emph{NO}(\varphi|_S)\longrightarrow0, \eqno(\star)$$
where $\mathscr{I}(X,S+B;\varphi)$ is the usual multiplier ideal sheaf associated to the weight $\varphi$ for the pair $(X,S+B)$ on $X$, and $i:S\hookrightarrow X,\ \iota$ together with $\rho$ are the natural inclusion and restriction map respectively.
\end{theorem}

\begin{remark}
$(1)$ By the construction of analytic adjoint ideal sheaf $\text{Adj}_S(X,S+B;\varphi)$, the surjectivity of $\rho$ relies heavily on the $L^2$ extension theorem and strong openness of multiplier ideals (cf. the proof of Theorem \ref{Aadjoint}). Moreover, we can deduce the coherence of $\text{Adj}_S(X,S+B;\varphi)$ from the adjunction exact sequence $(\star)$.

$(2)$ If $X$ is smooth, the above result has been established in \cite{Li_multiplier, Li_adjoint}; refer to \cite{Taka10, Eisen10} and \cite{Gue12, Kim15a, G-L_adjoint} for related topics.

$(3)$ Whenever $\varphi$ has analytic singularities, it follows that $\text{Adj}_S(X,S+B;\varphi)$ coincides with the algebraic adjoint ideal sheaf defined by Takagi and Eisenstein in \cite{Taka10, Eisen10}, which is independent of the choice of log resolution.
\end{remark}

\subsection{Inversion of adjunction}

In higher dimensional birational geometry, the adjunction/inversion of adjunction theorems have been important tools for studying singularities of pairs. Thanks to a crucial connectedness result (cf. \cite{Kollar92}, Theorem 17.4), the following inversion of adjunction was established by \cite{Shokurov92} in dimension 3 and by \cite{Kollar92} in general.

\begin{theorem} \emph{(\cite{Kollar92}, Theorem 17.6).} \label{Kollar_inversion}
Let $X$ be a normal complex space and $S\subseteq X$ an irreducible reduced Weil divisor. Let $B$ be an effective $\mathbb{Q}$-divisor on $X$ such that $K_X+S+B$ is $\mathbb{Q}$-Cartier. Then, $(X,S+B)$ is \emph{PLT} in a neighborhood of $S$ if and only if $(\widehat S,\emph{Diff}_{\widehat S}(B))$ is \emph{KLT}.
\end{theorem}

In a recent paper \cite{KimSeo21}, the authors present the following natural question posed by Koll\'ar whether we can find an analytic proof of Theorem \ref{Kollar_inversion} relying on an $L^2$ extension theorem. Concretely, we have

\begin{question} \emph{(\cite{KimSeo21}, Question B).} \label{Question_Kollar}
Does the inversion of adjunction Theorem \ref{Kollar_inversion} have an analytic proof using $L^2$ extension theorems?
\end{question}

As a corollary of Theorem \ref{Aadjoint}, we can deduce a generalized analytic inversion of adjunction as follows, by which we answer Question \ref{Question_Kollar} in the affirmative, i.e.,

\begin{theorem} \label{Answer_Kollar}
With the same hypotheses as before, then we have along $S$,
$$\emph{Adj}_S(X,S+B;\varphi)=\mathcal{O}_X\iff\mathscr{I}_\emph{NO}(\varphi|_S)=\mathcal{O}_S.\eqno{(\spadesuit)}$$
In particular, the solution to Question \ref{Question_Kollar} is positive.
\end{theorem}

As we will present, $(X,S+B)$ is PLT in a neighborhood of $S$ if and only if the analytic adjoint ideal sheaf $\text{Adj}_S(X,S+B;0)=\mathcal{O}_X$ along $S$, and $(\widehat S,\text{Diff}_{\widehat S}(B))$ is KLT if and only if the Nadel-Ohsawa multiplier ideal sheaf $\mathscr{I}_\text{NO}(0)=\mathcal{O}_S$. Therefore, in order to answer Question \ref{Question_Kollar}, it is equivalent to show that $\text{Adj}_S(X,S+B;0)=\mathcal{O}_X$ along $S$ if and only if $\mathscr{I}_\text{NO}(0)=\mathcal{O}_S$, which is an immediate consequence of $(\spadesuit)$ by taking a trivial weight $\varphi=0$.

\begin{remark}
When $X$ is smooth, Question \ref{Question_Kollar} has been answered implicitly by the adjunction exact sequence $(\star)$ established in \cite{Li_multiplier} or \cite{Li_adjoint}; see also Theorem 1.4 in \cite{KimSeo21} for an explicit presentation based on an analogous argument.
\end{remark}

\section{Preliminaries}

In this section, we will state some notions and useful results that will be used throughout this paper.

\begin{definition}
Let $X$ be a normal complex space and $\Delta=\sum d_k\Delta_k$ a Weil $\mathbb{R}$-divisor on $X$, where $\Delta_k$ are distinct prime divisors. The \emph{round-up} $\lceil\Delta\rceil$ and \emph{round-down} (or \emph{integral part}) $\lfloor\Delta\rfloor=[\Delta]$ of $\Delta$ are the integral divisors
$$\lceil\Delta\rceil=\sum\lceil d_k\rceil\Delta_k,\quad \lfloor\Delta\rfloor=\sum\lfloor d_k\rfloor\Delta_k,$$
where for any $c\in\mathbb{R}$ we denote the least integer $\geq c$ by $\lceil c\rceil$, and the largest integer $\leq c$ by $\lfloor c\rfloor=[c]$.

The \emph{fractional part} $\{\Delta\}$ and \emph{co-fractional part} $\langle\Delta\rangle$ of $\Delta$ are defined by $$\{\Delta\}=\Delta-\lfloor\Delta\rfloor,\quad \langle\Delta\rangle=\lceil\Delta\rceil-\Delta.$$
\end{definition}

\begin{definition} \label{MIS_Pair}
Let $(X,\Delta)$ be a log pair with $\varphi\in L_\text{loc}^1(X_\text{reg})$ . A positive definite \emph{adapted measure} $\upsilon$ on $X_\text{reg}$ for the pair $(X,\Delta)$ is defined as a positive measure locally of the form $\upsilon=e^{\phi}\cdot\Big((\sqrt{-1})^{\gamma n^2}\sigma\wedge\overline\sigma\Big)^{\frac{1}{\gamma}}$ on $X$, where $\phi$ is a bounded measurable function and $\sigma$ is a local generator of $\mathcal{O}_X\big(\gamma(K_X+\Delta)\big)$ on an open subset of $X$ for some real number $\gamma>0$ (cf. \cite{EGZ09}, Definition 6.5).

The \emph{multiplier ideal sheaf} associated to $\varphi$ on $X$ for the pair $(X,\Delta)$ is defined to be the $\mathcal{O}_X$-submudle $\mathscr{I}(X,\Delta;\varphi)\subseteq\mathscr{M}_X$ of germs of meromorphic functions $f\in\mathscr{M}_{X,x}$ such that $|f|^2e^{-2\varphi}$
is locally integrable at $x\in X$ with respect to an adapted measure $\upsilon$ defined as above. Here, $\varphi$ is regarded as the weight function and $\mathscr{I}(X,\Delta;\varphi)$ is independent of the choice of $\upsilon$; one can refer to \cite{La04} for an algebraic counterpart.

If $\varphi$ has analytic singularities on $X$, i.e., locally $\varphi=c\cdot\varphi_\mathfrak{a}:=c\cdot\log|\mathfrak{a}|+O(1)$ for some ideal sheaf $\mathfrak{a}\subseteq\mathcal{O}_X$ and $c\in\mathbb{R}$ on $X$, the multiplier ideal sheaf $\mathscr{I}(X,\Delta;\varphi)$ is also denoted by $\mathscr{I}(X,\Delta;c\cdot\varphi_\mathfrak{a})$ or $\mathscr{I}(X,\Delta;\mathfrak{a}^c)$. In the context, we associate every Weil $\mathbb{R}$-divisor $\Delta=\sum d_k\Delta_k$ with a natural weight $\varphi_\Delta$, which can be locally written as $\varphi_\Delta=\sum d_k\log|\mathscr{I}_{\Delta_k}|$.
\end{definition}

\begin{remark}
$(1)$ If $X$ is smooth and $\varphi\in\text{QPsh}(X)$, it follows that $\mathscr{I}(X,\Delta=0;\varphi)\subseteq\mathcal{O}_X$ is nothing but the multiplier ideal sheaf $\mathscr{I}(\varphi)$ on $X$ introduced by Nadel (cf. \cite{De10}).

$(2)$ If $\Delta$ is effective and $\varphi\in\text{QPsh}(X)$, then $\mathscr{I}(X,\Delta;\varphi)\subseteq\mathcal{O}_X$ is a coherent ideal sheaf.

$(3)$ If $\Delta$ is an $\mathbb{R}$-Cartier divisor on $X$, then we have $\mathscr{I}(X,\Delta;\varphi)=\mathscr{I}(X,0;\varphi+\varphi_\Delta)$.
\end{remark}

Similar to the smooth case, we have the following bimeromorphic transformation rule for pairs (see \cite{La04}, Proposition 9.3.62).

\begin{proposition}
Let $\pi:(\widetilde X,\Delta_{\widetilde X})\to(X,\Delta)$ be a bimeromorphic morphisms between log pairs such that $K_{\widetilde X}+\Delta_{\widetilde X}=\pi^*(K_X+\Delta)$ and $\pi_*\Delta_{\widetilde X}=\Delta$. Then, it follows that
$$\mathscr{I}(X,\Delta;\varphi)=\pi_*\left(\mathscr{I}(\widetilde X,\Delta_{\widetilde X};\varphi\circ\pi)\right).$$
Moreover, if $\widetilde X$ is smooth, we have
$$\mathscr{I}(X,\Delta;\varphi)=\pi_*\left(\mathcal{O}_{\widetilde X}(\lceil-\Delta_{\widetilde X}\rceil)\otimes\mathscr{I}(\widetilde X,0;\varphi\circ\pi+\varphi_{\langle-\Delta_{\widetilde X}\rangle})\right).$$
\end{proposition}

\begin{remark}
Let $(X,\Delta)$ be a log pair and the weight $\varphi_A=c_1\log|\mathfrak{a}_1|-c_2\log|\mathfrak{a}_2|$ be a difference of two quasi-psh functions with analytic singularities on $X$. Let $\pi:\widetilde X\to X$ be a common log resolution of the pair $(X,\Delta)$ and $\mathfrak{a}_k\ (k=1,2)$ with $\mathfrak{a}_k\cdot\mathcal{O}_{\widetilde X}=\mathcal{O}_{\widetilde X}(-D_k)$ and $K_{\widetilde X}-\pi^*(K_X+\Delta)+D_1+D_2$ having SNC supports. Then, it follows that
$$\mathscr{I}(X,\Delta;\varphi_A)=\pi_*\mathcal{O}_{\widetilde X}\left(\lceil K_{\widetilde X}-\pi^*(K_X+\Delta)-c_1D_1+c_2D_2\rceil\right).$$
\end{remark}

In the context, we will adopt a relative version of Grauert-Riemenschneider vanishing theorem for the higher direct images stated as below.

\begin{theorem} \emph{(\cite{Matsu_morphism}, Corollary 1.5).} \label{G-R_vanishing}
Let $\pi:X\to Y$ be a surjective proper (locally) K\"ahler morphism from a complex manifold $X$ to a complex space $Y$, and $(L,e^{-2\varphi_L})$ be a (possibly singular) Hermitian line bundle on $X$ with semi-positive curvature. Then, the higher direct image sheaf $$R^q\pi_*\Big(\omega_X\otimes\mathcal{O}_X(L)\otimes\mathscr{I}(\varphi_L)\Big)=0,$$
for every $q>\dim X-\dim Y$.
\end{theorem}

For the sake of convenience, we state a special case of the Ohsawa--Takegoshi-type $L^2$ extension theorem obtained by Guan and Zhou (Theorem 2.2 in \cite{G-Z_optimal}) as follows.

\begin{theorem} \label{OT}
Let $\Omega\subseteq\mathbb{C}^{n+r}$ be a pseudoconvex domain and $Z\subseteq\Omega$ an (closed) complex subspace of pure codimension $r$ such that $\mathscr{I}_Z$ is globally generated by holomorphic functions $g_1,\dots,g_m\in\mathcal{O}(\Omega)$.

Then, there is a constant $C>0$ (depending only on $r$) such that for any $\varphi\in\text{\emph{Psh}}(\Omega)$ and any $f\in\mathcal{O}_Z(Z)$ with $\int_Z|f|^2e^{-2\varphi}dV_{X}[\Psi]<+\infty$, there exists a holomorphic extension $\tilde f$ of $f$ to $\Omega$ such that
$$\int_\Omega\frac{|\tilde f|^2}{|g|^{2r}(-\log|g|)^2}e^{-2\varphi}dV_\Omega\leq C\cdot\int_Z|f|^2e^{-2\varphi}dV_{X}[\Psi],$$
where $\Psi=r\log|g|^2=r\log(|g_1|^2+\cdots+|g_m|^2)$.
\end{theorem}

In \cite{B-V03, Wlo08} and \cite{Bravo13, Eisen10}, the following desingularization theorem is established, which is necessary for our construction of the Nadel-Ohsawa multiplier/adjoint ideal sheaves for the higher codimensional case.

\begin{theorem} \emph{(Strong factorizing desingularization).} \label{SFD}
Let $X$ be a complex space of pure dimension and $Z\subseteq X$ a (closed) complex subspace with defining sheaf of ideals $\mathscr{I}_Z$ and no irreducible components contained in $X_\emph{sing}$. Then there exists an embedded resolution of singularities $\pi:\widetilde X\to X$ of $Z$ such that
$$\pi^*(\mathscr{I}_Z)=\mathscr{I}_{\widetilde Z}\cdot\mathscr{I}_{R_Z},$$
where $\mathscr{I}_{\widetilde Z}$ is the defining sheaf of ideals of the strict transform $\widetilde Z$ of $Z$ in $\widetilde X$ which has simple normal crossings with the exceptional divisor $\emph{Ex}(\pi)$ of $\pi$, and $\mathscr{I}_{R_Z}$ is the sheaf of ideals of a simple normal crossing divisor $R_Z$ supported on $\emph{Ex}(\pi)$.

Furthermore, we can choose $\pi$ to be a log resolution of any $\mathbb{R}$-linear combination of subschemes on $X$ not containing any component of $Z$ in its support (cf. \cite{Eisen10}, Corollary 3.2).
\end{theorem}

\section{Nadel-Ohsawa multiplier ideal sheaves on divisors for log pairs}

In this section, we present an Ohsawa-type measure on divisors arising from the research of so-called Ohsawa-Takigoshi $L^2$ extension theorem, by which will construct one type of multiplier ideal sheaves on divisors for log pairs.

\subsection{Ohsawa-type measure}

In order to establish a general $L^2$ extension theorem, Ohsawa \cite{Ohsawa5} introduced a positive measure on regular part of the associated complex subspace. Afterwards, associated with the same measure, the authors established various general $L^2$ extension theorems in \cite{G-Z_optimal, De16, CDM17} and so on. In the following, we will introduce an Ohsawa-type measure on divisors for log pairs to study singularities of pairs associated to psh functions in the sense of multiplier and adjoint ideal sheaves.

Let $(X,S+B)$ be a log pair with an adapted measure $\upsilon$ on $X_\text{reg}$, where $S\subseteq X$ is a reduced Weil divisor and $B\subseteq X$ is an $\mathbb{R}$-divisor on $X$ which has no irreducible components in common with $S$. Then, for any polar function $\Psi:X\to[-\infty,+\infty)$ along $S$ (i.e., locally written as $\log|\mathscr{I}_S|^2+O(1)$), we can associate a positive measure $dV_{S}[\Psi]$ on $S_\text{reg}$ defined by
$$\int_{S_\text{reg}}fdV_{S}[\Psi]=\limsup_{t\to+\infty}\int_{X_\text{reg}}f\cdot\mathbbm{1}_{\{-t-1<\Psi<-t\}}d\upsilon$$
for any nonnegative continuous function $f$ with $\text{Supp}\,f\,{\subset\subset}\,X_\text{reg}\backslash S_\text{sing}$, where $\mathbbm{1}_{\{-t-1<\Psi<-t\}}$ denotes the characteristic function of the set $\{z\in X\,|-t-1<\Psi(z)<-t\}$.

\begin{remark} \label{computation}
Let $\pi:\widetilde X\to X$ be a log resolution of $(X,S+B)$ or a strong factorizing desingularization of $S$ by Theorem \ref{SFD} such that $\pi^*(\mathscr{I}_S)=\mathscr{I}_ {\widetilde S}\cdot\mathscr{I}_{R_S}$, where $\widetilde S$ is the strict transform (may be disconnected) of $S$ in $\widetilde X$ and $R_S$ is an effective divisor supported on the exceptional divisor $\text{Ex}(\pi)$ of $\pi$. Write $$K_{\widetilde X}+\Delta_{\widetilde X}=\pi^*(K_X+S+B)=\pi^*(K_X+B)+\widetilde S+R_S.$$
Following the same argument as Remark 2.1 in \cite{Li_multiplier}, one can check that the measure $dV_{S}[\Psi]$ is locally the direct image of measures defined upstairs by
$$f\longmapsto\int_{\widetilde S}\left|f\circ\pi|_{\widetilde S}\right|^2\cdot\left|\frac{\text{Jac}(\pi)}{h_{R_S}}\Big|_{\widetilde S}\right|^2dV_{\widetilde S},$$
up to multiplicative bounded factors; where $\big(\text{Jac}(\pi)\big)=K_{\widetilde X}-\pi^*(K_X+B)=-\Delta_{\widetilde X}+\widetilde S+R_S$ and $h_{R_S}$ is the (local) defining function of $R_S$.
\end{remark}

\subsection{Nadel-Ohsawa multiplier ideals on divisors}

The usual analytic multiplier ideals are constructed by the integrability with respect to the Lebesgue measure or an adapted measure through the associated pluricanonical forms (e.g., Definition \ref{MIS_Pair}). In the analytic setting, the main difficulty of extending the notion of multiplier ideals to any singular complex (sub-)space is to choose a suitable measure in the sense of integrability. Thanks to the Ohsawa-type measure on divisors defined above, we are able to construct one type of multiplier ideals on divisors for log pairs in the following.

\begin{definition} \label{MIS_NO}
Let $(X,S+B)$ a log pair with an adapted measure $\upsilon$ on $X_\text{reg}$ as before. Let $\Psi:X\to[-\infty,+\infty)$ be a polar function along $S$ and  $\varphi\in L_\text{loc}^1(X_\text{reg})$ such that $\varphi|_S\not\equiv-\infty$ on every irreducible component of $S$.

Then, the \emph{Nadel-Ohsawa multiplier ideal sheaf} $\mathscr{I}_\text{NO}(\varphi|_S)$ on the divisor $S$ associated to weight $\varphi$ for the pair $(X,S+B)$ is defined to be the fractional ideal sheaf of germs of meromorphic functions $f\in\mathscr{M}_{S,x}$ such that $|f|^2e^{-2\varphi}$ is locally integrable at $x$ on $S$ with respect to the measure $dV_{S}[\Psi]$. One can check that $\mathscr{I}_\text{NO}(\varphi|_S)$ is independent of the choices of $\upsilon$ and $\Psi$.
\end{definition}

\begin{theorem} \label{property}
With the same notations as above and $\varphi\in\emph{QPsh}(X)$, then the Nadel-Ohsawa multiplier ideal sheaf $\mathscr{I}_\emph{NO}(\varphi|_S)\subseteq\mathscr{M}_{S}$ is a coherent fractional ideal sheaf and satisfies the strong openness property, i.e.,
$$\mathscr{I}_\emph{NO}(\varphi|_S)=\bigcup\limits_{\varepsilon>0}\mathscr{I}_\emph{NO}\big((1+\varepsilon)\varphi|_S\big).$$
\end{theorem}

\begin{proof}
Let $\pi:\widetilde X\to X$ be a log resolution of $(X,S+B)$ or a strong factorizing desingularization of $S$ by Theorem \ref{SFD} such that $\pi^*(\mathscr{I}_S)=\mathscr{I}_ {\widetilde S}\cdot\mathscr{I}_{R_S}$, where $\widetilde S$ is the strict transform (may be disconnected) of $X$ in $\widetilde X$ and $R_S$ is an effective divisor supported on the exceptional divisor $\text{Ex}(\pi)$ of $\pi$. Write $$K_{\widetilde X}+\Delta_{\widetilde X}=\pi^*(K_X+S+B)=\pi^*(K_X+B)+\widetilde S+R_S.$$
Then, by Remark \ref{computation} for any $f\in\mathscr{I}_\text{NO}(\varphi|_S)_x$ defined on a small enough neighborhood $U$ of $x$ in $X$, we have
\begin{equation*}
\begin{split}
+\infty>\int_U|f|^2e^{-2\varphi}dV_{X}[\Psi]=&\int_{\widetilde U=\pi|_{\widetilde S}^{-1}(U)}\left|f\circ\pi|_{\widetilde S}\right|^2\cdot e^{-2\varphi\circ\pi|_{\widetilde S}}\cdot\left|\frac{\text{Jac}(\pi)}{h_{R_S}}\Big|_{\widetilde S}\right|^2dV_{\widetilde S},
\end{split}
\end{equation*}
which implies that $$f\in\pi_*\left(\mathcal{O}_{\widetilde X}(\lceil-\Delta_{\widetilde X}+\widetilde S\rceil)|_{\widetilde S}\otimes\mathscr{I}(\varphi\circ\pi|_{\widetilde S}+\varphi_{\langle-\Delta_{\widetilde X}+\widetilde S\rangle}|_{\widetilde S})\right)_x,$$
where $\big(\text{Jac}(\pi)\big)=K_{\widetilde X}-\pi^*(K_X+B)=-\Delta_{\widetilde X}+\widetilde S+R_S$ and $h_{R_S}$ is the (local) defining function of $R_S$. On the other hand,
$$\pi_*\left(\mathcal{O}_{\widetilde X}(\lceil-\Delta_{\widetilde X}+\widetilde S\rceil)|_{\widetilde S}\otimes\mathscr{I}(\varphi\circ\pi|_{\widetilde S}+\varphi_{\langle-\Delta_{\widetilde X}+\widetilde S\rangle}|_{\widetilde S})\right)\subseteq\mathscr{M}_X$$
is an $\mathcal{O}_X$-submodule of germs of meromorphic functions $f\in\mathscr{M}_{S,x}$ such that $|f|^2e^{-2\varphi}$ is locally integrable at $x$ on $S$ with respect to the measure $dV_{S}[\Psi]$. Then, it follows that
$$\mathscr{I}_\text{NO}(\varphi|_S)=\pi_*\left(\mathcal{O}_{\widetilde X}(\lceil-\Delta_{\widetilde X}+\widetilde S\rceil)|_{\widetilde S}\otimes\mathscr{I}(\varphi\circ\pi|_{\widetilde S}+\varphi_{\langle-\Delta_{\widetilde X}+\widetilde S\rangle}|_{\widetilde S})\right),\eqno{(\heartsuit)}$$
which is coherent by Grauert's direct image theorem and Nadel's coherence theorem.
In addition, the strong openness of Nadel-Ohsawa multiplier ideals immediately follows from $(\heartsuit)$ and the result established by Guan and Zhou in \cite{G-Z_open} (see also \cite{Pham14}).
\end{proof}

\begin{remark}
$(1)$ When both of $X$ and $S$ are smooth, the Nadel-Ohsawa multiplier ideal sheaf $\mathscr{I}_\text{NO}(\varphi|_S)$ above is nothing but the usual multiplier ideal sheaf $\mathscr{I}(\varphi|_S+\varphi_B|_S)$ on $S$ introduced by Nadel.

$(2)$ In addition, if $\varphi_A$ is a quasi-psh function on $X$ possessing analytic singularities and $\varphi_A|_S\not\equiv-\infty$ on every irreducible component of $S$, by combining with an argument of log resolution, we can deduce that the fractional ideal sheaf $\mathscr{I}_\text{NO}\big((\varphi-\varphi_A)|_S\big)\subseteq\mathscr{M}_S$ is coherent.
\end{remark}

For the convenience of readers, we state the following result on Nadel-Ohsawa multiplier ideal sheaves established in \cite{Li_multiplier} for smooth ambient spaces (see also \cite{Li_adjoint} for the case of divisors), relying on the $L^2$ extension theorem and strong openness of multiplier ideal sheaves.

\begin{theorem} \emph{(\cite{Li_multiplier}, Theorem 1.4).} \label{Aadjoint_smooth}
Let $X$ be a complex manifold, $Z\subseteq X$ a complex subspace of pure codimension $r$ and $\varphi\in\emph{QPsh}(X)$ such that $\varphi|_Z\not\equiv-\infty$ on every irreducible component of $Z$.

Then, there exists an ideal sheaf $\emph{Adj}_Z(\varphi)\subseteq\mathcal{O}_X$ such that the following sequence of ideal sheaves is exact:
$$0\longrightarrow\mathscr{I}(\varphi+r\log|\mathscr{I}_Z|)\stackrel{\iota}{\longrightarrow}\emph{Adj}_Z(\varphi)\stackrel{\rho}{\longrightarrow} i_*\mathscr{I}_{\emph{NO}}(\varphi|_Z)\longrightarrow0.$$
\end{theorem}

\begin{remark} \label{Aadjoint_smoothRe}
In particular, when $r=1$ or $Z\subseteq X$ is smooth, the ideal sheaf $\text{Adj}_Z(\varphi)$ could be constructed as follows.

As the statement is local, we may assume that $X$ is a bounded Stein domain in some $\mathbb{C}^{N}$. Let $\mathscr{J}\subseteq\mathcal{O}_X$ be an ideal sheaf such that $\mathscr{J}|_Z=\mathscr{I}_\text{NO}(\varphi|_Z)\ $ (=$\mathscr{I}(\varphi|_Z)$ if $Z$ is smooth), which implies that $\mathscr{J}+\mathscr{I}_Z$ is independent of the choice of $\mathscr{J}$. Then, it follows that $$\text{Adj}_Z(\varphi):=\bigcup\limits_{\varepsilon>0}\text{Adj}_Z^0\Big((1+\varepsilon)\varphi\Big)\cap\Big(\mathscr{J}+\mathscr{I}_Z\Big),$$
where $\text{Adj}_Z^0(\varphi)\subseteq \mathcal{O}_X$ is an ideal sheaf of germs of holomorphic functions $f\in\mathcal{O}_{X,x}$ such that $$\frac{|f|^{2}e^{-2\varphi}}{|\mathscr{I}_Z|^{2r}\log^2|\mathscr{I}_Z|}$$ is locally integrable with respect to the Lebesgue measure near $x$ on $X$.
\end{remark}

\section{Proof of main results}

\subsection{Proof of Theroem \ref{Aadjoint}}

Let $\pi:\widetilde X\to X$ be a log resolution of $(X,S+B)$ or a strong factorizing desingularization of $S$ by Theorem \ref{SFD} such that $\pi^*(\mathscr{I}_S)=\mathscr{I}_ {\widetilde S}\cdot\mathscr{I}_{R_S}$, where $\widetilde S$ is the strict transform (may be disconnected) of $X$ in $\widetilde X$ and $R_S$ is an effective divisor supported on the exceptional divisor $\text{Ex}(\pi)$ of $\pi$. Write $$K_{\widetilde X}+\Delta_{\widetilde X}=\pi^*(K_X+S+B)=\pi^*(K_X+B)+\widetilde S+R_S.$$
Then, it follows from Theorem \ref{Aadjoint_smooth} that the following sequence of ideal sheaves is exact
\begin{equation*}
\begin{split}
0\longrightarrow\mathscr{I}(\varphi\circ\pi+\varphi_{\langle-\Delta_{\widetilde X}+\widetilde S\rangle})\otimes\mathcal{O}_{\widetilde X}(-\widetilde S)&\stackrel{\iota}{\longrightarrow}\text{Adj}_{\widetilde S}(\varphi\circ\pi+\varphi_{\langle-\Delta_{\widetilde X}+\widetilde S\rangle})\\
&\stackrel{\rho}{\longrightarrow} i_*\mathscr{I}(\varphi\circ\pi|_{\widetilde S}+\varphi_{\langle-\Delta_{\widetilde X}+\widetilde S\rangle}|_{\widetilde S})\longrightarrow0.
\end{split}
\end{equation*}
Here, the surjectivity of $\rho$ is due to the strong openness of multiplier ideals and an $L^2$ extension theorem, e.g., Theorem \ref{OT}.

Twist the exact sequence by $\mathcal{O}_{\widetilde X}(\lceil-\Delta_{\widetilde X}+\widetilde S\rceil)$, and then we deduce that
\begin{equation*}
\begin{split}
0\to\mathcal{O}_{\widetilde X}(\lceil-\Delta_{\widetilde X}\rceil)\otimes\mathscr{I}(\varphi\circ\pi+\varphi_{\langle-\Delta_{\widetilde X}+\widetilde S\rangle})&\to\text{Adj}_{\widetilde S}(\varphi\circ\pi+\varphi_{\langle-\Delta_{\widetilde X}+\widetilde S\rangle})\otimes\mathcal{O}_{\widetilde X}(\lceil-\Delta_{\widetilde X}+\widetilde S\rceil)\\
&\to i_*\mathscr{I}(\varphi\circ\pi|_{\widetilde S}+\varphi_{\langle-\Delta_{\widetilde X}+\widetilde S\rangle}|_{\widetilde S})\otimes\mathcal{O}_{\widetilde X}(\lceil-\Delta_{\widetilde X}+\widetilde S\rceil)\to0.
\end{split}
\end{equation*}
Taking
$$\text{Adj}_S(X,S+B;\varphi):=\pi_*\left(\text{Adj}_{\widetilde S}(\varphi\circ\pi+\varphi_{\langle-\Delta_{\widetilde X}+\widetilde S\rangle})\otimes\mathcal{O}_{\widetilde X}(\lceil-\Delta_{\widetilde X}+\widetilde S\rceil)\right)$$
and pushing forward the above exact sequence, we will infer the desired adjunction exact sequence
$$0\longrightarrow\mathscr{I}(X,S+B;\varphi)\stackrel{\iota}{\longrightarrow}\text{Adj}_S(X,S+B;\varphi)\stackrel{\rho}{\longrightarrow} i_*\mathscr{I}_\text{NO}(\varphi|_S)\longrightarrow0 \eqno(\star)$$
from the local vanishing (e.g., Theorem \ref{G-R_vanishing}) on the higher direct images of the term on the left; the proof is thereby concluded.
\hfill $\Box$

\subsection{Proof of Theorem \ref{Answer_Kollar}}

Note that $$\mathscr{I}(X,S+B;\varphi)\subseteq\mathscr{I}(X,S+B;0)\subseteq\mathscr{I}_S$$ by the definition. Then, by twisting the exact sequence $(\star)$ through by $\mathcal{O}_S$, we can deduce a restriction formula $\mathscr{I}_\text{NO}(\varphi|_S)=\text{Adj}_S(X,S+B;\varphi)\cdot\mathcal{O}_S$ from the right exactness of tensor functor. In particular, we have along $S$,
$$\text{Adj}_S(X,S+B;\varphi)=\mathcal{O}_X\iff\mathscr{I}_\text{NO}(\varphi|_S)=\mathcal{O}_S.\eqno{(\spadesuit)}$$

Let $\mathfrak{a}\subseteq\mathcal{O}_X$ be an ideal sheaf on $X$ with zeros not containing any irreducible component of $S$ and put $\varphi=c\cdot\varphi_{\mathfrak{a}}=c\cdot\log|\mathfrak{a}|$ for some $c\in\mathbb{R}_{\geq0}$. Let $\pi:\widetilde X\to X$ be a common log resolution of $(X,S+B)$ and $\mathfrak{a}$ such that $\pi^*(\mathscr{I}_S)=\mathscr{I}_ {\widetilde S}\cdot\mathscr{I}_{R_S}$ as before, and $\mathfrak{a}\cdot\mathcal{O}_{\widetilde X}=\mathcal{O}_{\widetilde X}(-F)$ for some effective divisor $F$ on $\widetilde X$. By a direct calculation of $\text{Adj}_{\widetilde S}(c\cdot\varphi_{\mathfrak{a}}\circ\pi+\varphi_{\langle-\Delta_{\widetilde X}+\widetilde S\rangle})$, we could obtain that $\text{Adj}_S(X,S+B;\mathfrak{a}^c)=\pi_*\mathcal{O}_{\widetilde X}\left(\lceil-\Delta_{\widetilde X}+\widetilde S-c\cdot F\rceil\right)$, which implies that $(X,S+B;\mathfrak{a}^c)$ is PLT in a neighborhood of $S$ if and only if the analytic adjoint ideal sheaf $\text{Adj}_S(X,S+B;\mathfrak{a}^c)=\mathcal{O}_X$ along $S$.

On the other hand, by the definition of different $\text{Diff}_{\widehat S}(B)$ it follows that $(\widehat S,\text{Diff}_{\widehat S}(B);\mathfrak{a}^c|_S)$ is KLT if and only if the divisor $\lceil-\Delta_{\widetilde X}+\widetilde S-c\cdot F\rceil\big|_{\widetilde S}$ is effective on $\widetilde S$, which is equivalent to triviality of the Nadel-Ohsawa multiplier ideal sheaf $\mathscr{I}_\text{NO}(\mathfrak{a}^c|_S)=\mathcal{O}_S$. Therefore, by taking $c=0$ it follows from $(\spadesuit)$ that $(X,S+B)$ is PLT in a neighborhood of $S$ if and only if $(\widehat S,\text{Diff}_{\widehat S}(B))$ is KLT.
\hfill $\Box$

\appendix
   \renewcommand{\appendixname}{Appendix~\Alph{section}}

\section{Analytic inversion of adjunction: higher codimension} \label{A}

In this appendix, analogous to the case of divisors, we introduce Nadel-Ohsawa multiplier and adjoint ideal sheaves along (closed) complex subspaces of higher codimension for log pairs, by which we establish an analytic inversion of adjunction in the higher codimensional case. When the ambient space is smooth, one can refer to \cite{Li_multiplier}.

Let $(X,\Delta)$ be a log pair with an adapted measure $\upsilon$ on $X_\text{reg}$ and $Z\subseteq X$ be a complex subspace of pure codimension $r$ such that there are no irreducible components of $Z$ contained in $X_\text{sing}\cup\text{Supp}\ \Delta$. As before, for any polar function $\Psi:X\to[-\infty,+\infty)$ along $S$ which can be locally written as $r\cdot\log|\mathscr{I}_Z|^2+O(1)$, we can associate a positive measure $dV_{Z}[\Psi]$ on $Z_\text{reg}$ defined by
$$\int_{Z_\text{reg}}fdV_{Z}[\Psi]=\limsup_{t\to+\infty}\int_{X_\text{reg}}fe^{-\Psi}\cdot\mathbbm{1}_{\{-t-1<\Psi<-t\}}d\upsilon$$
for any nonnegative continuous function $f$ with $\text{Supp}\,f\,{\subset\subset}\,X_\text{reg}\backslash Z_\text{sing}$, where $\mathbbm{1}_{\{-t-1<\Psi<-t\}}$ denotes the characteristic function of the set $\{z\in X\,|-t-1<\Psi(z)<-t\}$.

Let $\pi:\widetilde X\to X$ be a strong factorizing desingularization of $Z$ by Theorem \ref{SFD} such that $\pi^*(\mathscr{I}_Z)=\mathscr{I}_ {\widetilde Z}\cdot\mathscr{I}_{R_Z}$, where $\widetilde Z$ is the strict transform (may be disconnected) of $Z$ in $\widetilde X$ and $R_Z$ is an effective divisor supported on the exceptional divisor $\text{Ex}(\pi)$ of $\pi$. We in further obtain that the measure $dV_{Z}[\Psi]$ is locally the direct image of measures defined upstairs by
$$f\longmapsto\int_{\widetilde Z}\left|f\circ\pi|_{\widetilde Z}\right|^2\cdot\left|\frac{\text{Jac}(\pi)}{h_{R_Z}^r}\Big|_{\widetilde Z}\right|^2dV_{\widetilde Z},$$
up to multiplicative bounded factors; where $\big(\text{Jac}(\pi)\big)=K_{\widetilde X}-\pi^*(K_X+\Delta)=-\Delta_{\widetilde X}$ and $h_{R_Z}$ is the (local) defining function of $R_Z$.

\begin{definition}
With the same notations as above and $\varphi\in L_\text{loc}^1(X_\text{reg})$ such that $\varphi|_Z\not\equiv-\infty$ on every irreducible component of $Z$.
Then, the \emph{Nadel-Ohsawa multiplier ideal sheaf} $\mathscr{I}_\text{NO}(\varphi|_Z)$ on $Z$ associated to weight $\varphi$ for the pair $(X,\Delta)$ is defined to be the fractional ideal sheaf of germs of meromorphic functions $f\in\mathscr{M}_{Z,x}$ such that $|f|^2e^{-2\varphi}$ is locally integrable at $x$ on $Z$ with respect to the measure $dV_{Z}[\Psi]$. One can check that $\mathscr{I}_\text{NO}(\varphi|_Z)$ is independent of the choices of $\upsilon$ and $\Psi$.
\end{definition}

By the same argument as our proof of Theorem \ref{property}, we deduce that it also holds in the setting of higher codimension. Furthermore, we can establish the following analogue of Theorem \ref{Aadjoint} in the case:

\begin{theorem}
Let $(X,\Delta)$ be a log pair with $\Delta$ an effective $\mathbb{R}$-divisor, $Z\subseteq X$ a complex subspace of pure codimension $r$ such that there are no irreducible components of $Z$ contained in $X_\emph{sing}\cup\emph{Supp}\ \Delta$, and $\varphi\in\emph{QPsh}(X)$ satisfying $\varphi|_Z\not\equiv-\infty$ on every irreducible component of $Z$. Then, there exists an ideal sheaf
$$\emph{Adj}_Z(X,\Delta;\varphi)\subseteq\mathcal{O}_X,$$
called the \emph{analytic adjoint ideal sheaf} associated to the triple $(X,\Delta;\varphi)$ along $Z$, sitting in an exact sequence:
$$0\longrightarrow\mathscr{I}(X,\Delta;\varphi+r\log|\mathscr{I}_Z|)\stackrel{\iota}{\longrightarrow}\emph{Adj}_Z(X,\Delta;\varphi)\stackrel{\rho}{\longrightarrow} i_*\mathscr{I}_\emph{NO}(\varphi|_Z)\longrightarrow0,$$
where $i:Z\hookrightarrow X,\ \iota$ and $\rho$ are the natural inclusion and restriction map respectively.
\end{theorem}

\begin{proof}
One can make an argument by following the proof of Theorem \ref{Aadjoint} (or Theorem 1.4 in \cite{Li_multiplier}) line by line. On the construction of ideal sheaf $\text{Adj}_{\widetilde Z}(\psi)$ for any $\psi\in\text{QPsh}(\widetilde X)$, we refer to Remark \ref{Aadjoint_smoothRe}.
\end{proof}

\begin{remark}
The analytic adjoint ideal sheaf $\text{Adj}_Z(X,\Delta;\varphi)$ defined above coincides with the algebraic adjoint ideal sheaf defined by Takagi and Eisenstein in \cite{Taka10, Eisen10} whenever $\varphi$ has analytic singularities.
\end{remark}

As an immediate consequence similar to the case of divisors, we obtain the following analytic inversion of adjunction for higher codimension.

\begin{theorem}
With the same hypotheses as above, then it follows that along $Z$,
$$\emph{Adj}_Z(X,\Delta;\varphi)=\mathcal{O}_X\iff\mathscr{I}_\emph{NO}(\varphi|_Z)=\mathcal{O}_Z.$$
\end{theorem}

\end{document}